\documentclass[11pt]{amsart}
\usepackage{amsmath}
\usepackage{amssymb}
\usepackage{euscript}
\usepackage{pstricks}
\usepackage{eepic}
\usepackage{times}
\usepackage[all]{xy}
\numberwithin{equation}{section}
\usepackage[font=small,labelfont=bf]{caption}
\numberwithin{figure}{section}

\newtheorem{definition}[equation]{Definition}
\newtheorem{theorem}[equation]{Theorem}
\newtheorem{proposition}[equation]{Proposition}
\newtheorem{corollary}[equation]{Corollary}

\newtheorem{remark}[equation]{Remark}
\newtheorem{conjecture}[equation]{Conjecture}

\newcommand{\w}{\wedge}
\newcommand{\n}{\notag}
\newcommand{\noi}{\noindent}

\newcommand{\hook}{\hookrightarrow}

\newcommand{\one}{\vspace{1mm}}
\newcommand{\two}{\vspace{2mm}}
\newcommand{\sq}{ \;  \square}

\newcommand{\C}{\mathbb{C}}

\newcommand{\Q}{\mathbb{Q}}

\newcommand{\W}{\mathcal{W}}

\newcommand{\PP}{\mathbb{P}}

\newcommand{\x } {\textnormal{x}}

\newcommand{\del}{\delta \phi}

\newcommand{\Pf} {\noi\emph{Proof.}$\; \;$}

\newcommand{\beq}{\begin{equation}}
\newcommand{\eeq}{\end{equation} }

\begin{document}

\sloppy

\title{Legendrian Gronwall conjecture}
\author{Joe S. Wang}
\address{Danyang, Corea 395-806}
\email{jswang12@gmail.com}
\keywords{Legendrian 3-web, linearization, Gronwall conjecture}
\subjclass[2000]{53A60}

\begin{abstract}
{
The Gronwall conjecture states that a planar 3-web of foliations
which admits more than one distinct linearizations is locally equivalent to an algebraic web.
We propose an analogue of the Gronwall conjecture
for the 3-web of foliations by Legendrian curves in a contact three manifold.
The Legendrian Gronwall conjecture states that a Legendrian 3-web admits
at most one distinct local linearization, with the only exception
when it is locally equivalent to the dual linear Legendrian 3-web
of the Legendrian twisted cubic in $\,\PP^3$.
We give a partial answer to the conjecture in the affirmative
for the class of Legendrian 3-webs of maximum rank.
We also show that a linear Legendrian 3-web which is sufficiently flat at a reference point
is rigid under local linear Legendrian deformation.
}
\end{abstract}

\maketitle

\section{Introduction}\label{sec1}
Let $\, M$ be a connected contact three manifold.
A Legendrian $\, d$-web on $\, M$ is by definition
a set of $\, d$ pairwise transversal foliations of $\, M$ by Legendrian curves.
The Legendrian web was introduced in \cite{Wa2} for a second order generalization
of the classical planar web.

Abelian relations and rank, two of the central concepts in web geometry,
are analogously defined for the Legendrian web.
The main result of \cite{Wa2} was that
the rank of a Legendrian $\,d$-web admits the optimal bound
$\, \rho_d=\frac{(d-1)(d-2)(2d+3)}{6}.$
We also gave an analytic characterization of the Legendrian 3-webs of maximum rank three.

\begin{figure}[htp]
\centering
\begin{picture}(100, 60)(117, 2)

\put(161,40){$\textnormal{Z}$}
\put(175,25){$\searrow$}
\put(145,25){$\swarrow$}
\put(184,10){ $\mathbb{Q}^3$ }
\put(129,10){ $\, \PP^3$ }
\put(185,30){$\pi_1$}
\put(135,30){$\pi_0$}
\end{picture}

\caption{Projective duality  associated with  \,Sp$_2 \C$}
\label{11duality}
\end{figure}
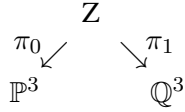

In this paper we study the linearization problem for the Legendrian 3-webs
with a more algebro-geometric perspective.
As described in \cite{Wa2}, the algebraic model for the Legendrian web theory
is provided by the projective duality associated with the simple Lie group Sp$_2\C$,
Figure \ref{11duality}, see Section \ref{sec11}, \cite{Br}.
By the standard dual construction, a null degree $\, d$ surface $\, \Sigma\subset\Q^3$
induces the dual $\, d$-web of Legendrian lines on a generic small
open subset of $\,\PP^3$.\footnotemark
\footnotetext{An analytic surface $\, \Sigma\subset\Q^3\subset\PP^4$ has null degree $\, d$
when it intersects a generic null line of $\, \Q^3$ at $\, d$ points.
Hence it has degree $\, 2d$ as a surface in $\,\PP^4$.}
Generalizing this, a \emph{linearization} of a Legendrian $\, d$-web
on a contact three manifold $\, M$ is defined as a contactomorphism $\, M \hook \PP^3$
for which the image of each leaf of the Legendrian foliations is mapped to a Legendrian line.
A natural question arises as to which Legendrian $\, d$-webs are linearizable,
and how unique such a linearization is.
The linearizability problem for a planar web can be traced back to Balschke,
\cite{BB,AGL,GL,PP,Pir} and the references therein.
See \cite{Wa} for the references on the Gronwall conjecture.

It is the uniqueness part of the linearization problem that we are interested in.
In planar web geometry, the Gronwall conjecture states that
a planar 3-web admits more than one distinct linearizations (uniqueness fails)
whenever the 3-web is locally equivalent to an algebraic web.
In other words, the conjecture claims that the failure of unique linearization,
or equivalently the linear deformability, implies that the planar 3-web is
essentially algebraic.

One of the motivations for the present work was therefore the idea
that the condition for linear deformability may lead to a special, possibly algebraic,
class of linear Legendrian 3-webs. Through these examples one might hope to gain an insight
on the analogues of Abel's theorem and its converse in our setting, Section \ref{sec4}.

\two
For the Legendrian $\, d$-webs, $\, d \geq 4$,
the uniqueness of linearization follows from the local normal form
for the projectively flat third order ODE's, \cite[Corollary 4.2]{SY}.
The case of Legendrian 3-webs on the other hand is not obvious.
We shall employ the method of moving frames and the local differential analysis
to analyze the linear Legendrian deformation and rigidity of
the linear Legendrian 3-webs on an open subset of $\,\PP^3$.

\two
\textbf{Main results}.

1.
The number of distinct local linearizations of a Legendrian $\, 3$-web
is uniformly bounded, Theorem \ref{31thm} (the bound is far from being optimal).

2.
A linear Legendrian 3-web on a connected open subset of $\,\PP^3$ is rigid under
linear Legendrian deformation when it is sufficiently flat at a single point,
Theorem \ref{32thm}.

3.
We propose the Legendrian analogue of the Gronwall conjecture, Conjecture \ref{conj}.
The conjecture claims that a Legendrian 3-web on a connected contact three manifold
admits at most one distinct linearization, with the only exception when
the 3-web is locally equivalent to the dual linear Legendrian 3-web of the Legendrian
twisted cubic curve in $\,\PP^3$
(in which case it has exactly two distinct local linearizations).
We verify the conjecture for the class of Legendrian 3-webs of maximum rank three,
Theorem \ref{4thm}.

\two
Let us give a description of the deformable Legendrian 3-web in \textbf{Main results} 3.
Let $\,\gamma \subset \PP^3$ be the Legendrian twisted cubic. Take a generic point
$\,\x \in \PP^3.$ Let $\,\PP^2_{\x}$ be the contact 2-plane at $\,\x$,
which is the union of the Legendrian lines through $\,\x$.
Since $\,\gamma$ is a curve of degree 3, $\,\PP^2_{\x}$ intersects $\,\gamma$ at
3 points. By definition of $\,\PP^2_{\x}$
these points determine the corresponding 3\, Legendrian lines
through $\,\x$. It is clear that as $\,\x$ varies this construction defines
a linear Legendrian $\,3$-web $\,\W_{\gamma}$ on a generic small open subset of $\, \PP^3$.
Our analysis shows that $\,\W_{\gamma}$ admits another (exactly one more)
distinct linearization, and that it is the unique such linear Legendrian 3-web
of maximum rank.

The three Abelian relations of $\,\W_{\gamma}$ admit the following geometrical
interpretation as a generalized addition law for the Legendrian twisted cubic.
Take a set of 3 generic points $\,\textnormal{p}_i, \, i=1,\,2,\,3,$ on $\,\gamma.$
Then the 3 planes $\,\PP^2_{\textnormal{p}_i}$ intersect at a single point,
and one consequently gets a set of 3 concurrent Legendrian lines.
The choice of three points on a given curve $\,\gamma$ depends on three one-dimensional
parameters. Once these points are fixed, the choice of a Legendrian line through
each of these points also depends on three one-dimensional parameters.\footnotemark
\footnotetext{Here we mean a choice of a section of a $\,\PP^1\times\PP^1\times\PP^1$-bundle
over $\,\gamma\times\gamma\times\gamma.$}
It is evident from the above incidence relation that, roughly speaking,
the latter three parameters are functions of the former three parameters.
The three Abelian relations of $\, \W_{\gamma}$ implies that in fact
there exist three additive form of functional relations among the total
six one-dimensional parameters.

Considering the obvious analogue of the converse of Abel's theorem, Section \ref{sec4},
the relevant point here would be that the incidence geometry described above
manifests itself through the \emph{additive} functional equations
among the related geometric quantities
which by construction are the first integrals of the foliations of $\,\W_{\gamma}$.
As is the case with the planar web geometry,
perhaps this additivity provides a concrete evidence to infer,
and in turn it would imply via the converse of Abel's theorem,
the underlying algebraic structure.

\one
There is another algebro-geometric implication of the three Abelian relations of
$\,\W_{\gamma}$ on the degeneration of K3 surfaces in $\,\PP^4$
in the context of the presumed Legendrian analogue of Abel's theorem and its converse.
Consider the surface $\,\Sigma_{\gamma}=\pi_1\circ\pi_0^{-1}(\gamma)\subset\Q^3$.
Since $\,\gamma$ is Legendrian, $\,\Sigma_{\gamma}$ is the tangent developable
of the Gauss map of $\,\gamma$, which is a null rational normal curve in $\,\Q^3$.
Via the converse of Legendrian Abel's theorem, the three Abelian relations of $\,\W_{\gamma}$
would then imply that $\,\Sigma_{\gamma}$ supports three generalized closed
holomorphic 1-forms. On the other hand, $\,\Sigma_{\gamma}$
lies in the intersection of $\,\Q^3$ with a cubic hypersurface of $\,\PP^4$.
A smooth complete intersection of type $\,(2,\,3)$ in $\,\PP^4$ is a K3 surface,
which has no nonzero holomorphic 1-forms. The analytic surface $\,\Sigma_{\gamma}$
represents in this way a degeneration of K3 surfaces where
the dimension of the space of closed holomorphic 1-forms jumps by three.

\one
Let us give an outline of the paper.
In Section \ref{sec11}, we record the generalities on
the duality between the projective space $\,\PP^3$ and the 3-quadric $\,\Q^3\subset\PP^4$.
In Section \ref{sec2}, the moving frame method is applied to
determine a normalized frame bundle associated with a linear Legendrian 3-web
on an open subset of $\,\PP^3$. The normalized frame bundle serves as a base point
for the linear Legendrian deformation of the 3-web.
In Section \ref{sec3}, the standard deformation analysis in terms of the deformed Maurer-Cartan
equation leads to the closed structure equation for the deformation parameters.
The differential compatibility conditions among these parameters imply
a set of polynomial equations that they must satisfy. The rigidity results
Theorem \ref{31thm} and Theorem \ref{32thm} are the immediate consequences
of the analysis of the root structure of the first few
low degree polynomial compatibility equations.
In Section \ref{sec4}, we implement the procedure established in Section \ref{sec3}
to the class of maximum rank Legendrian 3-webs.
Drawing from \cite{Wa2} the analytic characterization of such 3-webs,
we first observe that a Legendrian 3-web of maximum rank always admits
a linearization as the dual 3-web of a union of three hyperplane sections in $\,\Q^3$.
A further analysis of the polynomial compatibility equations shows that
there exists a unique linear Legendrian 3-web of this kind
which is contact equivalent (and not projectively equivalent)
to the dual 3-web of the Legendrian twisted cubic described above,
and that otherwise these 3-webs are rigid under
linear Legendrian deformation. Based on this analysis,
we propose the Legendrian analogue of the converse of Abel's theorem.
In Section \ref{sec5},
we state the Legendrian Gronwall conjecture,
and present a few remarks on the related problems.

\one
We assume the complex analytic category.
The analysis and results are valid within the real smooth category
with only minor modifications.
The moving frame method, and the over-determined PDE machinery
are used throughout the paper.
We refer the reader to \cite{BCG3}\cite{IvLa} for the standard references.

\subsection{Projective duality}\label{sec11}
Let $\, V = \C^{4}$ be the four dimensional complex vector space.
Let $\, \varpi$ be the standard symplectic 2-form on $\, V$.
Let $\,\PP^3 = \PP(V)$ be the projectivization equipped with the induced
contact structure. The contact 2-plane field $\, \mathcal{H}$ on
$\, \PP^3$ is defined by
\beq\label{11Hdefi}
\mathcal{H}_{\mbox{x}} =
\PP((\hat{\mbox{x}} \lrcorner \, \varpi )^{\perp}), \; \;
\mbox{for}\;\; \mbox{x} \in \PP^3,
\eeq
where   $ \hat{\mbox{x}} \in  V$ is any de-projectivization of  x.
Since $\,\varpi$ is non-degenerate,
$( \hat{\mbox{x}}  \lrcorner \,\varpi )^{\perp} \subset V$ is a codimension one subspace
containing $ \hat{\mbox{x}}$, and its projectivization
$\, \PP((\hat{\mbox{x}} \lrcorner \, \varpi )^{\perp}) \subset \PP^3$
is a hyperplane at x. The symplectic group Sp$_2 \C$ acts transitively on $\, \PP^3$
as a  group of contact transformation.

Let $\, \bigwedge_0^2V = \langle\,\varpi\,\rangle^{\perp} \subset  \bigwedge^2V$
be the five dimensional subspace of isotropic 2-vectors.
Let $\, Lag(V)$ be the set  of  two dimensional Lagrangian subspaces of $\, V$.
$\, Lag(V)$ is identified with the 3-quadric $\, \Q^3 \subset \PP^4=\PP(\bigwedge_0^2V)$
via the Pl\"ucker embedding.
Sp$_2 \C$ acts transitively on $\, \Q^3$, and
$\,\Q^3$ inherits the Sp$_2\C$ invariant non-degenerate conformal structure.

Let $\, \textnormal{Z}=\PP(\mathcal{H}) \to \PP^3$ be the bundle of
Legendrian(\,$\mathcal{H}$-horizontal) line elements.
It can be defined as the incidence space
$$
\textnormal{Z}=\{(\x,\,\textnormal{L}) \in\PP^3\times\Q^3\;\,|\;\;
\hat{\x}\w\hat{\textnormal{L}}=0\;\,\}.
$$
The projective duality in Figure \ref{11duality} represents the following
Sp$_2\C$ equivariant incidence double fibration, Figure \ref{11double}.
By definition of Z, it is clear that both fibers of $\, \pi_0, \, \pi_1$ are
isomorphic to $\, \PP^1$.
A fiber of $\, \pi_0$, the set of Legendrian lines through a point in $\,\PP^3$,
projects to a null line of $\, \Q^3$,
and dually a fiber of $\, \pi_1$, the set of null lines through a point in $\,\Q^3$,
projects to a Legendrian line of $\, \PP^3$.
The duality correspondence can be summarized by;

\one
\emph{\qquad
$\, \PP^3$ is the space of null lines in $\, \Q^3$,
and dually $\,\Q^3$ is the space of Legendrian lines in $\,\PP^3$.}

\one
There exists an immediate application of the duality principle.
Let $\, \gamma \subset \PP^3$ be a Legendrian ($\, \mathcal{H}$-horizontal) curve.
The Gauss map $\, \hat{\gamma}\subset\Q^3$ is defined as
the tangent map of $\,\gamma$. By duality,
$\, \hat{\gamma}$ is the envelope of a one parameter family of null lines.
Hence $\, \hat{\gamma}$ itself is a null curve in $\,\Q^3.$

\begin{figure}[htp]
\centering
\begin{picture}(300,89)(20,0)
\put(150,72){$\mbox{Sp$_2\C$}$} \put(172,57){$\pi$}
\put(162,57){$\downarrow$}
\put(161,40){$\textnormal{Z} = \mbox{Sp$_2\C$}/ \textnormal{P}$}
\put(175,25){$\searrow$} \put(145,25){$\swarrow$}
\put(184,10){ $\Q^3=Lag(V)$ } \put(129,10){ $ \PP^3$ }
\put(187,27){$\pi_1$} \put(133,27){$\pi_0$}
\end{picture}
\caption{Incidence double fibration}
\label{11double}
\end{figure}

\one
To fix the notation, let us define the projection
maps $\, \pi,  \, \pi_0$, and  $\, \pi_1 \, $ explicitly.
Let  $\,Z = (Z_0, \, Z_1, \, Z_2, \, Z_3)$ denote the Sp$_2\C \subset $ SL$_4\C$\, frame
of $\, V$ such that  the 2-vector
$\, \varpi_{\flat} = Z_0 \w Z_2 + Z_1 \w Z_3$ is dual to the symplectic form $\, \varpi$.
Define
\begin{align}\label{11defipi}
\pi (Z) &= ([Z_0], \, [Z_0 \w Z_1]),   \\
\pi_0 ([Z_0], \, [Z_0 \w Z_1]) &= [Z_0],  \n \\
\pi_1 ([Z_0], \, [Z_0 \w Z_1]) &= [Z_0 \w Z_1]. \n
\end{align}
\noi In this formulation, the stabilizer subgroup $\, \textnormal{P}$ in Figure \ref{11double}
is of the form
\beq \label{11group}
\textnormal{P}  = \{ \,
\begin{pmatrix}
A       & B \\
\cdot  & (A^{t})^{-1}
\end{pmatrix} \},
\eeq where $\, (A^{-1} B) ^t = A^{-1} B$, and \beq A  = \{ \,
\begin{pmatrix}
*       & *  \\
\cdot & *
\end{pmatrix} \}. \n
\eeq Here '$\cdot$' denotes 0  and '$*$' is arbitrary.

The Sp$_2\C$-frame $\, Z$  satisfies the structure equation
\beq\label{11Zframe}
d  Z = Z \, \phi
\eeq
for the Maurer-Cartan form $\, \phi$  of  Sp$_2\C$.
$\,\phi$ satisfies the structure equation
\beq\label{11MC}
d \phi+\phi \w \phi=0.
\eeq
The components of $\, \phi$ are  denoted by
\beq\label{11strt}
\phi =
\begin{pmatrix}
\alpha &  \gamma \\
\beta   &  - \alpha^t
\end{pmatrix}
=
\begin{pmatrix}
\alpha^0_0 & \alpha^0_1   & \gamma^0_0  & \gamma^0_1  \\
\alpha^1_0 & \alpha^1_1   & \gamma^1_0  & \gamma^1_1  \\
\beta^0_0 & \beta^0_1     & -\alpha^0_0 & -\alpha^1_0 \\
\beta^1_0 & \beta^1_1     &-\alpha^0_1  & -\alpha^1_1
\end{pmatrix},
\eeq
where $\, \{ \, \alpha, \, \beta, \, \gamma \, \}$ are 2-by-2
matrix 1-forms such that $\, \beta^t=\beta, \, \gamma^t=\gamma$.

\section{Linear Legendrian 3-webs}\label{sec2}
Let $\, \W$ be a Legendrian $\,d$-web on a contact three manifold $\, M$, \cite{Wa2}.
A linearization of $\, \W$ is a contactomorphism $\, M \hook \PP^3$
such that each leaf of the foliations is mapped to a Legendrian line.
Two linearizations are equivalent when they are isomorphic up to
projective transformation by Sp$_2\C$, and otherwise distinct. More specifically;

\begin{definition}\label{1webdefi}
Let $\, \PP^3$ be equipped with the standard \,\textnormal{Sp}$_2\C$ invariant
homogeneous contact structure, Section \ref{sec11}.
Let $\, M \subset \PP^3$ be a connected open subset.
A \emph{linear Legendrian $\,d$-web} on $\, M$
is a set of $\,d$ pairwise transversal foliations of $\, M$
by Legendrian lines.
\end{definition}

There exists a distinguished class of linear Legendrian webs.
Recall that an analytic surface $\, \Sigma \subset \Q^3$ has
null degree $\, d$ when it intersects a generic null line of $\,\Q^3$ at $\, d$ points.
\begin{definition}\label{2algwebdefi}
Let $\, \W$ be a linear Legendrian $\, d$-web on an open subset of $\,\PP^3$.
$\, \W$ is \emph{algebraic}, and is denoted by $\, \W_{\Sigma}$,
when it is induced from a null degree $\, d$
analytic surface $\, \Sigma \subset \Q^3$ by the standard dual construction.
\end{definition}

Let us give a description of a particular class of algebraic Legendrian webs
induced from a curve in $\,\PP^3$.
Let $\,\gamma\subset\PP^3$ be a degree $\,d$-curve.
Take a generic point $\,\x\in\PP^3$. The contact 2-plane $\,\mathcal{H}_{\x}$
intersects $\,\gamma$ at $\,d$ points $\, \textnormal{p}_i(\x),\,i=1,\,2,\,...\,d$.
By definition of $\,\mathcal{H}_{\x}$, \eqref{11Hdefi},
one gets a set of $\, d$ Legendrian lines through $\,\x$.
It is clear then that this incidence construction defines a linear Legendrian
$\,d$-web $\, \W_{\gamma}$ on a generic small open subset of $\, \PP^3$.
Note from Figure \ref{11double} that
$\,\W_{\gamma}=\W_{\pi_1\circ\,\pi_0^{-1}(\gamma)}.$

The dual Legendrian web $\,\W_{\gamma}$ obtained in this way inherits a set of Abelian relations
from the holomorphic 1-forms on $\,\gamma$.
Let $\, \Omega\in H^0(\gamma,\,\Omega^1)$ be a holomorphic 1-form.
By Hartogs' theorem, the trace of $\,\Omega$,
$$
\textnormal{Tr}\,\Omega = \sum\textnormal{p}_i^*\Omega, \n
$$
 is a holomorphic 1-form on $\,\PP^3$,
and hence must vanish identically.\footnotemark
\footnotetext{The trace is obtained by the pull back of
the map $\,\PP^3\smallsetminus\gamma \to \gamma^{(d)}$.}
Each 1-form $\,\textnormal{p}_i^*\Omega$ trivially vanishes on
the Legendrian line determined by $\,\textnormal{p}_i$.
This implies that there exists a linear map
from $\, H^0(\gamma,\,\Omega^1)$ to the space of
Abelian relations of $\, \W_{\gamma}$.

Consider the special case when $\, \gamma$ is itself a Legendrian curve.
The associated surface $\,\pi_1\circ\pi_0^{-1}(\gamma)$ becomes the tangent developable
of the Gauss map of $\,\gamma$.
We shall see that in this case there are Abelian relations of $\,\W_{\gamma}$
which do not come from $\, H^0(\gamma,\,\Omega^1)$, Section \ref{sec4}.

\subsection{Equivalence problem}\label{sec21}
Let $\, \W$ be a linear Legendrian 3-web.
For definiteness,
we assume the three Legendrian foliations of $\, \W$ are ordered,
and denote them by
$\, \W  = \cup_{i=1}^3 \mathcal{F}^i$.
Note by duality Figure \ref{11duality} that
each foliation $\, \mathcal{F}^i$
corresponds to a (possibly singular) surface $\, \Sigma^i \subset \mathbb{Q}^3$.
An immersed surface in the three manifold  $\, \Q^3$ is locally described
as a graph of one scalar function of two variables.
One may argue that
the local moduli space of linear Legendrian 3-webs in $\, \PP^3$
depends on three arbitrary scalar functions of two variables.

In this section,
we apply the method of moving frames
and determine the $\, \textnormal{Sp}_2\C$ invariant structure equation
for a linear Legendrian 3-web.
The analysis will result in a principal bundle $\, \mathcal{B} \to M$ equipped with
a normalized $\, \mathfrak{sp}_2\C$-valued Maurer-Cartan form $\, \phi$.
The functional relations among the coefficients of the components of $\, \phi$
are the basic local invariants of a linear Legendrian 3-web.
Since the data $\, ( \, \mathcal{B}, \, \phi \,)$ are canonically associated
with the given linear Legendrian 3-web,
they will serve as a reference point for the problem of
deformation and rigidity to be discussed in Section \ref{sec3}.

\one
Let $\, \W$ be a linear Legendrian 3-web on an open subset $\, M \subset \PP^3$.
Let $\,\textnormal{B} \subset \textnormal{Sp}_2\C \to M$ denote the induced
principal right $\,\textnormal{P}$-bundle.
We continue the analysis from Section \ref{sec11}.

\one
\textbf{Step 0}.
From \eqref{11strt}, set $\, \alpha^1_0 = \omega^1, \,\beta^1_0 = \omega^2;
     \, \beta^0_0 = 2 \theta; \, \alpha^0_0 = \rho_0$.
The Maurer-Cartan form $\, \phi$ is written as
\beq\label{2strt0}
\phi=
\begin{pmatrix}
\rho_0 & \alpha^0_1   & \gamma^0_0  & \gamma^0_1  \\
\omega^1   & \alpha^1_1   & \gamma^1_0  & \gamma^1_1  \\
2 \theta   & \omega^2     & -\rho_0 & -\omega^1 \\
\omega^2   & \beta^1_1    &-\alpha^0_1  & -\alpha^1_1
\end{pmatrix}.
\eeq
From the general theory of moving frames,
one may apply the fiber group action by $\, \textnormal{P} \subset \textnormal{Sp}_2\C$
to arrange so that the three linear Legendrian foliations are defined by
\beq\label{2defining}
\mathcal{F}^i = \langle\, \omega^i, \, \theta \, \, \rangle^{\perp}, \; i = 1, \, 2, \, 3,
\eeq
where $\, \omega^3 = -(\omega^1 + \omega^2)$.

Let $\,\mathcal{B}_0 \subset \textnormal{B}  \to M$
be the sub-bundle defined by \eqref{2defining}.
Assuming that the foliations are ordered,
the structure group $\, P_0 \subset \textnormal{P}$ of $\, \mathcal{B}_0$
is reduced to
$$
P_0 = \{
\begin{pmatrix}
*      & *     & *  & * \\
\cdot  & \pm  1 & *  & \cdot \\
\cdot& \cdot  & * & \cdot  \\
\cdot& \cdot  & * & \pm  1
\end{pmatrix}   \},
$$
where '$\, \cdot$ '  denotes  0.
On $\, \mathcal{B}_0 \to M$,
the  1-forms $\, \beta^1_1, \, \gamma^1_1; \, \alpha^1_1$
are semi-basic, and one may write
\begin{align}\label{2semibasic}
 \beta^1_1 &=\epsilon_1 \omega^1+\epsilon_2 \omega^2+ A_{1} \theta,   \\
 \gamma^1_1&=\epsilon_3 \omega^1+\epsilon_4 \omega^2+ A_{2} \theta, \n \\
 \alpha^1_1&=\epsilon_5 \omega^1+\epsilon_6 \omega^2+ A_{9} \theta, \n
\end{align}
for the coefficients $\, \epsilon_j; \, A_{k}.$

\one
\textbf{Step 1}.
The condition that $\, \W$ is linear imposes a set of relations on $\, \epsilon_j$'s.
By  \eqref{11Zframe}, \eqref{2defining}, one must have
\begin{align}
dZ_0, \, dZ_1 & \equiv 0, \mod \; \; Z_0, \, Z_1; \quad \omega^2, \, \theta, \n \\
dZ_0, \, dZ_3 & \equiv 0, \mod \; \; Z_0, \, Z_3; \quad \omega^1, \, \theta, \n \\
dZ_0, \, d(Z_1-Z_3) & \equiv 0, \mod \; \; Z_0, \, Z_1-Z_3; \quad \omega^3, \, \theta. \n
\end{align}
A computation shows that this implies in \eqref{2semibasic}
\beq\label{2linearity}
\epsilon_1 = 0, \; \; \epsilon_4 = 0,  \; \;
2(\epsilon_5 - \epsilon_6) - \epsilon_2 - \epsilon_3=0.
\eeq

\one
\textbf{Step 2}.
$\, \phi$ satisfies the  structure equation $\,d\phi+\phi \w \phi=0$.
For a notational purpose, denote $\, \Phi = d\phi+\phi \w \phi$,
which must vanish identically.
For instance,
$\, \Phi^1_0, \, \Phi^3_0; \, \Phi^2_0\,$ give the formulae for the exterior derivatives
$\, d\omega^1, \, d\omega^2; \, d\theta$, etc.

$\, \Phi^1_3\w \omega^2\w \theta, \, \Phi^3_1 \w \omega^1\w \theta$ show that
$$
\begin{array}{lll}
d\epsilon_2 &\equiv 2 \alpha^0_1, &\mod \;\;  \omega^1, \, \omega^2, \, \theta, \, \rho_0, \n \\
d\epsilon_3 &\equiv 2 \gamma^0_1, &\mod \;\;  \omega^1, \, \omega^2, \, \theta, \, \rho_0. \n
\end{array}
$$
One may apply the fiber group action by $\, P_0$ to translate
so that
\beq\label{2B1}
 \epsilon_2, \, \epsilon_3 = 0.
\eeq
From \eqref{2linearity}, denote
$\, \epsilon_5 = \epsilon_6 = A_0.$
The exterior derivative of $\, A_0$ is written as
$$dA_0 = -A_0 \rho_0 + A_{0,1} \omega^1 +  A_{0,2} \omega^2 + A_{0,0} \theta.$$
We shall adopt the similar notation for the covariant derivative of a coefficient
for the rest of the paper.

Let $\, \mathcal{B}_1 \subset \mathcal{B}_0$ be the sub-bundle defined by \eqref{2B1}.
The structure group $\, P_1 \subset P_0$ of $\, \mathcal{B}_1$
is reduced to the form
$$
P_1 = \{
\begin{pmatrix}
*      & \cdot     & *  & \cdot \\
\cdot  & \pm  1 & \cdot  & \cdot \\
\cdot& \cdot  & * & \cdot  \\
\cdot& \cdot  & \cdot & \pm  1
\end{pmatrix}   \}.
$$
On $\, \mathcal{B}_1 \to M$,
the  1-forms $\, \alpha^0_1, \, \gamma^0_1$
are semi-basic.
$\, \Phi^1_3 \w \theta, \, \Phi^3_1 \w \theta$ show that
one may write
\begin{align}\label{2semibasic2}
 \alpha^0_1 &=\frac{A_{1}}{2} \omega^1  +\epsilon_7 \omega^2  + B_{1} \theta,   \\
 \gamma^0_1 &=\epsilon_8 \omega^1 - \frac{A_{2}}{2} \omega^2  + B_{2} \theta, \n
\end{align}
for the coefficients $\, \epsilon_7, \, \epsilon_8; \, B_{1},\, B_{2}$.

$\, \Phi^1_1 \w \theta$ with this relation  gives
$\,A_{9}=-\epsilon_7-\epsilon_8-A_{0,1}+A_{0,2}-2 A_0^2.$

\one
\textbf{Step 3}.
$\, \omega^1 \w \Phi^0_1 + \Phi^0_3 \w \omega^2$ gives
$$d(\epsilon_7-\epsilon_8)  \equiv - 2 \gamma^0_0,
\mod \; \; \omega^1, \, \omega^2, \, \theta, \, \rho_0.$$
One may translate, and denote
\beq\label{2B2}
\epsilon_7 = \epsilon_8 = B_0.
\eeq

Let $\, \mathcal{B} \subset \mathcal{B}_1$ be the sub-bundle defined by \eqref{2B2}.
The structure group $\, P \subset P_1$ of $\, \mathcal{B}$
is reduced to
\beq\label{2groupP2}
P = \{
\begin{pmatrix}
*      & \cdot     & \cdot  & \cdot \\
\cdot  & \pm  1 & \cdot  & \cdot \\
\cdot& \cdot  & * & \cdot  \\
\cdot& \cdot  & \cdot & \pm  1
\end{pmatrix}   \}.
\eeq
On $\, \mathcal{B} \to M$,
the  1-form  $\,  \gamma^0_0$  is  semi-basic.
One may write
\begin{align}\label{2semibasic3}
 \gamma^0_0 &= C_1 \omega^1 + C_2 \omega^2+ C_9 \theta,
\end{align}
for the coefficients $\, C_1, \, C_2, \, C_9$.\footnotemark
\footnotetext{
We use the subscript '9' in $\, C_9$ in place of '0'
to indicate that $\, C_9$ has the higher scaling weight than $\, C_1, \, C_2$
under the action of the structure group $\, P$,
see \eqref{2structureeq}.}

\one
\textbf{Step 4}.
Differentiating \eqref{2semibasic}, \eqref{2semibasic2}, \eqref{2semibasic3},
and
examining the rest of the components of $\, \Phi$,
one obtains the following structure equations.\footnotemark
\footnotetext{
The method of differential analysis used here is referred to as
the \emph{prolongation}, \cite{BCG3}.
It is the process of successively adding the derivatives as the new variables,
under the condition of contact which indicates that these new variables are the derivatives.
It allows one to access the differential relations (not necessarily of higher order)
which are possibly hidden and can only be detected by
examining the higher order derivatives.
The computation was carried out using the \textbf{Maple}.
}
\beq\label{2structureeq}
\small{
\begin{array}{rl}
dA_0&=-A_{{0}}\rho_0 +A_{{0,1}}\omega^1+A_{{0,2}}\omega^2+A_{0,0}\theta,   \\
dA_{0,1}&\equiv -2\,A_{{0,1}}\rho_0 + A_{{0,1,1}}\omega^1
+\left( A_{{0,1,1}}+6\,A_{{0}}A_{{0,1}}+4\,A_{{0}}B_{{0}}
+2\, A_{{0}}^{3}+2\,A_{{0,0}}-2\,A_{{1}}A_{{0}}-5\,B_{{1}}-2\,C_{{1}} \right) \omega^2,   \\
dA_{0,2}&\equiv -2\,A_{{0,2}}\rho_0
+\left( A_{{0,1,1}}-2\,C_{{1}}+4\,A_{{0}}B_{{0}}+5\,A_{{0}}A_{{0,1}}
-A_{{0}}A_{{0,2}}-2\,A_{{1}}A_{{0}}+2\, A_{{0}}^{3}+A_{{0,0}}-5\,B_{{1}} \right) \omega^1   \\
&\quad+ \left(A_{{0,1,1}}+ 5\,B_{{2}}+2\,C_{{2}}+5\,A_{{0}}A_{{0,2}}-2\,A_{{0}}A_{{2}}
-2\,C_{{1}}-5\,B_{{1}}+3\,A_{{0,0}}-2\,A_{{1}}A_{{0}}
+5\,A_{{0}}A_{{0,1}} \right) \omega^2,   \\
dA_{0,0}&\equiv -3\,A_{{0,0}}\rho_0
+  \left( A_{{0,1,0}}+A_{{0}}B_{{1}}-2\,A_{{0}}C_{{1}}-2\,A_{{0,2}}A_{{1}}
+4\,A_{{0,1}}B_{{0}}-A_{{0,1}}A_{{0,2}}+2\,A_{{0,1}}A_{{0}}^{2}+A_{{0,1}}^{2}\right)\omega^1   \\
&\quad + \left( A_{{0,2,0}}-4\,A_{{0,2}}B_{{0}}-2\,A_{{0,1}}A_{{2}}+A_{{0}}B_{{2}}
-2\,A_{{0}}C_{{2}}-A_{{0,1}}A_{{0,2}}
-2\,A_{{0,2}}A_{{0}}^{2}+A_{{0,2}}^{2} \right) \omega^2,   \\
dA_1&\equiv -2\,A_{{1}}\rho_0 + 2 A_{{1}}A_{{0}}\omega^1
+ \left( -2\,B_{{1}}+2\,A_{{1}}A_{{0}} \right) \omega^2,   \\
dA_2&\equiv -2\,A_{{2}}\rho_0 + \left( -2\,A_{{0}}A_{{2}}-2\,B_{{2}} \right) \omega^1
-2 A_{{0} }A_{{2}}\omega^2,   \\
dB_0&\equiv -2\,B_{{0}}\rho_0 + \left( -C_{{1}}-2\,B_{{1}} \right) \omega^1
+ \left( C_{{2}}+2\,B_{{2}} \right) \omega^2,   \\
dB_1&\equiv -3\,B_{{1}}\rho_0 +\left( -A_{{0,2}}A_{{1}}+A_{{1}}A_{{0,1}}
+2\,A_{{1}}A_{{0}}^{2}+A_{{0}}B_{{1}}+\frac{1}{2}\,A_{{1,0}} \right) \omega^1   \\
&\quad+ \left( -4\,B_{{0}}^{2}+2\,C_{{9}}+A_{{0}}B_{{2}}+B_{{2,1}}
-A_{{1}}A_{{2}}+A_{{0}}B_{{1}} \right) \omega^2,   \\
dB_2&\equiv-3\,B_{{2}}\rho_0+B_{{2,1}}\omega^{{1}}
+ \left( 2\,A_{{0}}^{2}A_{{2}}-A_{{2}}A_{{0,2}}-A_{{0}}B_{{2}}
-\frac{1}{2}\,A_{{2,0}}+A_{{0,1}}A_{{2}} \right)\omega^{{2}},  \\
dC_1&\equiv -3\,C_{{1}}\rho_0 +C_{{1,1}}\omega^1
+ \left( -\frac{5}{2}\,B_{{2,1}}+A_{{0}}C_{{1}}+A_{{1}}A_{{2}}-\frac{5}{2}\,A_{{0}}B_{{2}}
-2\,C_{{9}}+4\,B_0^{2} \right) \omega^2,   \\
dC_2&\equiv -3\,C_{{2}}\rho_0
+ \left( -\frac{5}{2}\,B_{{2,1}}-A_{{0}}C_{{2}}-\frac{5}{2}\,A_{{0}}B_{{2}}
-3\,C_{{9}}+6\,B_{{0}}^{2}+\frac{3}{2}\,A_{{1}}A_{{2}} \right) \omega^1
+C_{{2,2}}\omega^2,   \\
dC_9&\equiv -4\,C_{{9}}\rho_0 +\left( -B_{{2}}A_{{1}}+2\,B_{{1}}B_{{0}}+C_{{1,0}}
+4\,C_{{1}}B_{{0}}+C_{{1}}A_{{0,1}}-C_{{1}}A_{{0,2}}+2\,C_{{1}}A_{{0}}^{2}
-2\,C_{{2}}A_{{1}} \right) \omega^1   \\
&\quad + \left( -2\,B_{{2}}B_{{0}}-B_{{1}}A_{{2}}+C_{{2,0}}-2\,C_{{1}}A_{{2}}
-4\,C_{{2}}B_{{0}}-A_{{0,1}}C_{{2}}+C_{{2}}A_{{0,2}}
-2\,C_{{2}}A_{{0}}^{2} \right) \omega^2,  \\
dB_{2,1}&\equiv -4\,B_{{2,1}}\rho_0
+ ( -2\,B_{{1,0}}-4\,B_{{1}}A_{{0}}^{2}+2\,B_{{1}}A_{{0,2}}-2\,C_{{1,0}}
-4\,C_{{1}}A_{{0}}^{2}+2\,C_{{1}}A_{{0,2}}-2\,B_{{1}}A_{{0,1}}-2\,C_{{1}}A_{{0,1}}  \\
&\quad-16\,C_{{1}}B_{{0}}-16\,B_{{1}}B_{{0}}-4\,B_{{2}}A_{{1}}-A_{{0,1}}B_{{2}}
-A_{{0}}B_{{2,1}}  ) \omega^1
+  ( B_{{2}}A_{{0,2}}-2\,A_{{0,1}}B_{{2}}+\frac{1}{2}\,A_{{0}}A_{{2,0}}
-3\,A_{{0}}^{2}B_{{2}}   \\
&\quad+4\,B_{{2}}B_{{0}}+2\,B_{{2,0}}-2\,A_{{0}}^{3}A_{{2}}
+4\,C_{{1}}A_{{2}}
+4\,B_{{1}}A_{{2}}-A_{{0}}A_{{0,1}}A_{{2}}
+A_{{0}}A_{{0,2}}A_{{2}}  ) \omega^2, \mod \; \theta.
\end{array}}
\eeq

\normalsize
\noindent
The $\,\theta$-derivative terms, e.g.,
$\, A_{0,1,0}, A_{1,0}, B_{1,0}, \,... \,$, are all independent
with the one exception that
$$
B_{0,0}=C_{{9}}-2\, B_{{0}}^{2}+A_{{0}}B_{{2}}+B_{{2,1}}-\frac{1}{2}\,A_{{1}}A_{{2}}.
$$
\begin{remark}
In the language of the  theory of differential systems,
this set of structure equations is involutive
and a general analytic solution
(linear Legendrian 3-web with the given structure equations)
depends on three arbitrary functions of two variables as expected,
see the remark at the beginning of this section.
See \cite{BCG3} for the details.
\end{remark}
\begin{proposition}\label{2propo}
Let $\, \W$ be a linear Legendrian 3-web
on a connected open subset $\, M \subset \PP^3$.
There exists a canonically associated principal bundle
$\, \mathcal{B} \subset \textnormal{Sp}_2\C \to M$
with the structure group \eqref{2groupP2}.
The induced $\, \mathfrak{sp}_2\C$-valued Maurer-Cartan form $\, \phi$, \eqref{2strt0},
is normalized on $\, \mathcal{B}$ such that
\beq\label{2phi}
\phi= \left[ \begin {array}{cccc}
 \rho_0&\frac{1}{2} A_{{1}}\omega^1+B_{{0}}\omega^2
 +B_{{1}}\theta&C_{{1}}\omega^1+C_{{2}}\omega^2+C_{{9}}\theta
 &B_{{0}}\omega^1-\frac{1}{2}A_{{2}}\omega^2+B_{{2}}\theta\\
\noalign{\medskip}\omega^1&\phi^1_1
&B_{{0}}\omega^1-\frac{1}{2}A_{{2}}\omega^2+B_{{2}}\theta&A_{{2}}\theta\\
\noalign{\medskip}2\,\theta&\omega^2&-\rho_0&-\omega^1\\
\noalign{\medskip}\omega^2&A_{{1}}\theta&-\frac{1}{2}A_{{1}}\omega^1
-B_{{0}}\omega^2-B_{{1}}\theta&-\phi^1_1
\end {array} \right],
\eeq
where
$\, \phi^1_1= -A_{{0}}\omega^3
+(A_{{0,2}}-2\, B_{{0}}- A_{{0,1}}-2  A_{{0}}^{2})\theta$.
The structure coefficients $\, A_i, \, B_j, \, C_k$ and their derivatives
satisfy  \eqref{2structureeq}.

Two linear Legendrian 3-webs $\, \W, \, \W'$ are congruent up to
\textnormal{Sp}$_2\C$ motion
whenever the corresponding data $\, (\mathcal{B}, \, \phi)$ and
 $\, (\mathcal{B}', \, \phi')$ are isomorphic.
\end{proposition}

Let us rephrase the argument at the beginning of this section
with a view to applying Proposition \ref{2propo}.
Given a linear Legendrian 3-web $\, \W$ on $\, M \subset \PP^3$,
it determines 3 sections $\, M^i \subset \textnormal{Z}, \, i = 1, \, 2, \, 3,$
by definition of the duality in Figure \ref{11double}.
The linearity of $\, \W$ implies that each $\, M^i$ is tangent to
the fibers of $\, \pi_1$. Under the projection by $\, \pi_1$,  $\, M^i$ is mapped to
a surface $\, \Sigma^i \subset \Q^3$.
The image $\, \pi_1\circ\pi_0^{-1}(\x)$ is the dual null line of $\, \x \in M$
which intersects $\, \Sigma = \cup_{i=1}^3 \Sigma^i$ at 3 points.
The local geometry of a linear Legendrian 3-web in this way corresponds to
the semi-global geometry of a union of 3 pieces of surfaces in $\, \Q^3$.
Before we proceed to the problem of deformation,
let us consider an example where this dual interpretation allows
a simple description of a class of  linear Legendrian 3-webs.

Suppose for a linear Legendrian 3-web
the relative invariant $\, A_0$ vanishes identically,
$$A_0 \equiv 0.$$
An analysis shows that in this case the Maurer-Cartan form $\, \phi\,$
reduces to
\beq\label{2linearQ}
\phi =
 \left[ \begin {array}{cccc} \rho_0&\frac{1}{2} A_{{1}}\omega^1+B_{{0}}\omega^2
 & (2 B_{{0}}^{2} +\frac{1}{2}A_{{1}}A_{{2}}) \theta &B_{{0}}\omega^1
-\frac{1}{2} A_{{2}}\omega^2\\\noalign{\medskip}\omega^1&-2\,B_{{0}}\theta
&B_{{0}}\omega^1-\frac{1}{2} A_{{2}}\omega^2&A_{{2}}\theta\\\noalign{\medskip}2\,
\theta&\omega^2&-\rho_0&-\omega^1\\\noalign{\medskip}\omega^2&A_{{1}}\theta
&-\frac{1}{2} A_{{1}}\omega^1-B_{{0}}\omega^2&2 B_{{0}}\theta\end {array}
 \right],
\eeq
where $\, dA_1 = -2A_1\rho_0, \,  dA_2 = -2A_2\rho_0, \,  dB_0 = -2B_0\rho_0$.

Choose the Legendrian foliation $\, \mathcal{F}^1$ defined by
$\, \langle \, \omega^2, \, \theta \, \rangle^{\perp}$.
The corresponding surface $\, \Sigma^1 \subset \Q^3 \subset \PP(\bigwedge^2_0\C^4)$
is described by $\, [ Z_0 \w Z_1 ]$ (here we follow the notation
from Section \ref{sec11}).
A direct computation by successively differentiating $\, [ Z_0 \w Z_1 ]$
shows that $\, \Sigma^1$ is a part of a hyperplane section $\, H^1 \subset \Q^3$.
From the similar analysis for the foliations $\, \mathcal{F}^2, \, \mathcal{F}^3$,
one concludes that;

\one
\emph{Let $\, \W$ be a linear Legendrian 3-web with vanishing relative invariant $\, A_0$.
Then $\, \W$ is a part of the algebraic Legendrian 3-web $\, \W_{\Sigma}$
induced by $\, \Sigma = \cup_{i=1}^3 H^i$, a union of \,3 hyperplane sections in $\, \Q^3$.}

\one
We shall see in Section \ref{sec4}
that this class of Legendrian 3-webs account for all the linear Legendrian 3-webs
of maximum rank, with the only exception of the dual 3-web
of the Legendrian twisted cubic curve.

\section{Deformation, and rigidity}\label{sec3}
In this section, we establish the fundamental structure equation
for the linear Legendrian deformation of a linear Legendrian 3-web.
A variant of the moving frame method is applied,
and the analysis leads to the closed structure equation for the three deformation parameters.
The differential compatibility conditions of this structure equation
generate a sequence of polynomial equations that the deformation parameters must satisfy.

We currently have a partial understanding of the root structure of these polynomial equations.
An elementary examination of the first few polynomials shows that;
1) the number of distinct linearizations of a Legendrian $\,3$-web is uniformly bounded,
Theorem \ref{31thm},\,
2) if the Legendrian 3-web is sufficiently flat at a point, it admits at most one
distinct local linearization, Theorem \ref{32thm}.

\one
Let $\,\W$ be a linear Legendrian 3-web on a connected open subset $\,M\subset\PP^3$.
Let $\,\mathcal{B}\to M$ be the associated adapted bundle
with the $\, \mathfrak{sp}_2\C$-valued normalized Maurer-Cartan form $\, \phi$,
Proposition \ref{2propo}.
Let $\,M\hook\PP^3$ be another distinct linearization of $\,\W$.
Let $\,\textnormal{B}'\to M$ be the associated P$'\simeq\textnormal{P}$-bundle with
the Maurer-Cartan form $\, \phi'$. We employ the method of moving frames to normalize
the frame bundle $\,\textnormal{B}'$ based at $\,\mathcal{B}$.
The structure equation for the difference $\,\phi'-\phi $ then gives
the aforementioned polynomial compatibility equations.

Set \beq\label{3setpi}
\phi' = \phi + \del.
\eeq
The components of $\, \del$  are denoted by
$$
\del =
\begin{pmatrix}
\delta \alpha & \delta\gamma \\
\delta\beta  & - \delta \alpha^t
\end{pmatrix}.
$$
For the notational purpose, set
\beq\label{3fundast}
\Delta = d(\del) + \del \w \phi + \phi \w \del + \del \w \del.
\eeq
Maurer-Cartan equations  for $\, \phi'$ and $\, \phi$ imply that
$\,\Delta$ must vanish identically.

\one
\textbf{Step 0$'$}.
Applying the fiber group action by
$\,\textnormal{P}'\subset \textnormal{Sp}_2\C$
as in \textbf{Step 0} of Section \ref{sec2},
one may translate
\beq\label{3step00}
 \delta\alpha^1_0 = 0, \,\delta\beta^1_0 = 0;\, \delta\beta^0_0 = 0.
\eeq
Under these relations, the structure group is reduced to
$$
\textnormal{P}'_0 = \{ \pm
\begin{pmatrix}
1      & \cdot     & *  & \cdot \\
\cdot  & 1   & \cdot  & \cdot \\
\cdot& \cdot  & 1 & \cdot  \\
\cdot& \cdot  & \cdot & 1
\end{pmatrix}   \}.
$$

$\Delta^2_0$  gives $\delta\alpha^0_0 \w \theta=0.$
Applying the fiber group action by $\, \textnormal{P}'_0$,
one may translate
\beq\label{3step01}
 \delta\alpha^0_0 = 0.
\eeq
Under this relation, the structure group is reduced to
the center of Sp$_2\C$, $Z(\textnormal{Sp}_2\C) = \{ \, \pm I_4 \, \}$.
The $\mathfrak{sp}_2\C$-valued deformation 1-form $\, \delta \phi$ becomes
\beq
\delta \phi=
\begin{pmatrix}
\cdot   & *   & *  & *  \\
\cdot   & *   & *  & *  \\
\cdot   & \cdot        & \cdot & \cdot \\
\cdot   & *    & *  &*
\end{pmatrix}.\n
\eeq
Starting from this initial state,
by successive applications of \eqref{3fundast},
we intend to find the compatibility equations for the linear Legendrian deformation.

Differentiating \eqref{3step00}, \eqref{3step01},
one gets
$$
\begin{pmatrix}
 \delta \alpha^1_1 & \delta \gamma^1_1 & 2\delta\gamma^0_1 \\
-\delta  \beta^1_1 & \delta \alpha^1_1 & 2\delta\alpha^0_1 \\
 \delta \alpha^0_1 & \delta \gamma^0_1 & 2\delta\gamma^0_0
\end{pmatrix}\w
\begin{pmatrix}
\omega^1 \\ \omega^2 \\ \theta
\end{pmatrix}=0.
$$
By the Cartan's lemma, one may write
\beq\label{3step0prep}
\begin{pmatrix}
 \delta \alpha^1_1 \\ \delta \gamma^1_1 \\ \delta \beta^1_1 \\
 \delta\alpha^0_1  \\ \delta\gamma^0_1 \\ \delta\gamma^0_0
\end{pmatrix} =
\begin{pmatrix}
\mu_1  & \mu_2 & 2\mu_9   \\  \mu_2   & \mu_4 &  2\mu_{10}   \\
\mu_5  &-\mu_1 &-2\mu_7   \\  \mu_7   & \mu_9 &  2\mu_{11}    \\
\mu_9  & \mu_{10} & 2\mu_{12} \\  \mu_{11}  & \mu_{12} &  \mu_{11,0}
\end{pmatrix}
\begin{pmatrix}
\omega^1 \\ \omega^2 \\ \theta
\end{pmatrix},\n
\eeq
for the coefficients $\, \mu_l$.

\one
\textbf{Step 1$'$}.
The condition that it is a linear Legendrian deformation
imposes a set of relations on $\, \mu_k$'s.
A computation similar as in \textbf{Step 1} of Section \ref{sec2} shows that
$\,
\mu_4 = 0, \, \mu_5 = 0,  \, \mu_1-\mu_2=0.
$
We set accordingly
\beq\label{3step0set}
\begin{pmatrix}
 \delta \alpha^1_1 \\ \delta \gamma^1_1 \\ \delta \beta^1_1 \\
 \delta\alpha^0_1  \\ \delta\gamma^0_1 \\ \delta\gamma^0_0
\end{pmatrix} =
\begin{pmatrix}
a_0  & a_0  & 2b_0   \\   a_0   & \cdot  &  2b_3   \\
\cdot   &- a_0 & -2b_1  \\  b_1   & b_0 &  2 c_1   \\
b_0  & b_3 & 2 c_2 \\  c_1  &  c_2  & c_9
\end{pmatrix}
\begin{pmatrix}
\omega^1 \\ \omega^2 \\ \theta
\end{pmatrix}.
\eeq
The derivatives of these coefficients will be denoted similarly as before, e.g.,
$  da_0 = -a_0\rho_0+a_{0,1}\omega^1 + a_{0,2}\omega^2 + a_{0,0}\theta.$

\one

\textbf{Step 2$'$}.
$\, \Delta^1_3\w   \theta, \, \Delta^3_1  \w \theta,  \,\Delta^1_1  \w \theta$
give
\beq\label{3step0set0}
\begin{array}{ll}
a_{0,1}&=-4\,b_{{1}}+a_{{0}}A_{{0}}+2\,a_{{0}}^{2},   \\
a_{0,2}&=\;\;\;4\,b_{{3}}-a_{{0}}A_{{0}}-2\,a_{{0}}^{2},   \\
b_0&=b_1+b_3-a_{{0}}A_{{0}}- \frac{3}{4} \,a_{{0}}^{2}.
\end{array}\eeq

The remaining equations from
$\, \Delta^1_3, \, \Delta^3_1,  \,\Delta^1_1$
give
\beq\label{3step0set1}
\begin{array}{ll}
b_{1,1}&=2\,A_{{0}}b_{{1}}-2\,A_{{1}}a_{{0}}+2\,b_{{1}}a_{{0}},   \\
b_{1,2}&=2\,c_{{1}}+2\,A_{{0}}a_{{0}}^{2}+\frac{3}{2}\, a_{{0}}^{3}-2\,a_{{0}}b_{{3}}
-\frac{1}{2}\,A_{{0,2}}a_{{0}}+2\,A_{{0}}b_{{1}}+\frac{1}{2} \,A_{{0,1}}a_{{0}}
-A_{{1}}a_{{0}}+A_{{0}}^{2}a_{{0}}+\frac{1}{2} \,a_{{0,0}},   \\
b_{3,1}&=-2\,c_{{2}}-2\,A_{{0}}a_{{0}}^{2}-\frac{3}{2}\,a_{{0}}^{3}
+2\,b_{{1}}a_{{0}}-A_{{0}}^{2}a_{{0}}-2\,A_{{0}}b_{{3}}
-A_{{2}}a_{{0}}-\frac{1}{2}\,A_{{0,1}}a_{{0}}+\frac{1}{2}\,A_{{0,2}}a_{{0}}
+\frac{1}{2} \,a_{{0,0}},   \\
b_{3,2}&=-2\,A_{{0}}b_{{3}}-2\,A_{{2}}a_{{0}}-2\,b_3 a_{{0}},   \\
c_{1}&=-\frac{9}{2}\,a_{{0}}^{3}-\frac{11}{2}\,A_{{0}} a_{{0}}^{2}+3\,b_{{1}}a_{{0}}
+\frac{2}{3}\,A_{{0}}b_{{1}}+\frac{10}{3}\,A_{{0}}b_{{3}}+6\,a_{{0}}b_{{3}}
-3\, A_{{0}}^{2}a_{{0}}-\frac{4}{3}\,A_{{0,1}}a_{{0}}   \\
&\quad +\frac{5}{3}\,A_{{0,2}}a_{{0}}-2\,B_{{0}}a_{{0}}+\frac{1}{2}\,A_{{1}}a_{{0}},   \\
c_{2}&=c_1-3\,a_{{0}}b_{{3}}-\frac{8}{3}\,A_{{0}}b_{{3}}-\frac{1}{2}\,A_{{2}}a_{{0}}
-\frac{1}{3}\,A_{{0,1}}a_{{0}}-\frac{1}{3}\,A_{{0,2}}a_{{0}}+3\,b_{{1}}a_{{0}}
+\frac{8}{3}\,A_{{0}}b_{{1}}-\frac{1}{2}\,A_{{1}}a_{{0}}.
\end{array}\eeq

The compatibility equations from
$\,d(d(a_0))=0$
give
\beq\label{3step0set2}
\begin{array}{ll}
a_{0,0}&=-\frac{1}{3}\,A_{{0,1}}a_{{0}}-\frac{1}{3}\,A_{{0,2}}a_{{0}}
+\frac{8}{3}\,A_{{0}}b_{{1}}-\frac{8}{3}\,A_{{0}}b_{{3}},   \\
b_{1,0}&=5\,A_{{0}} a_{{0}}^{3}+ \left( -\frac{2}{3} \,A_{{0,2}}
+7\, A_{{0}}^{2}+A_{{1}} +\frac{1}{3} \,A_{{0,1}}+2\,B_{{0}} \right)  a_{{0}}^{2}
+ \mathcal{O}(b_1a_0, b_3a_0; \, a_0, b_1, b_3),  \\
b_{3,0}&=-5\,A_{{0}} a_{{0}}^{3}+ \left( -\frac{2}{3}\,A_{{0,1}}-7\, A_{{0}}^{2}
+A_{{2}}+\frac{1}{3}\,A_{{0,2}} -2\,B_{{0}} \right) a_{{0}}^{2}
+ \mathcal{O}(b_1a_0, b_3a_0; \, a_0, b_1, b_3).
\end{array}\eeq
Here $\, \mathcal{O}(b_1a_0, b_3a_0; \, a_0, b_1, b_3)$ means the terms that are
linear combination of $\, \{ \, b_1a_0, b_3a_0; \, a_0, b_1, b_3 \, \}$
with the coefficients in $\, A_i, \ B_j, \, C_k$'s and their derivatives.

At this stage, there are three components $\, \Delta^0_1, \, \Delta^0_3; \, \Delta^0_2$
left to be checked.

\one

\textbf{Step 3$'$}.
$\,\Delta^0_1\w\omega^1$ finally gives
\beq\label{3step0set3}
\begin{array}{ll}
c_9 &=
{\frac {81}{8}}\, a_{{0}}^{4}+19\,A_{{0}}a_{{0}}^{3}
+ \left( -26\,A_{{2}}+27\,b_{{1}}-4\,A_{{1}}-105\,b_{{3}}-\frac{5}{2}\,A_{{0,2}}
-19\,B_{{0}}-\frac{19}{2}\,A_{{0,1}}+5\,A_{{0}}^{2} \right)a_{{0}}^{2}   \\
&\quad +2\,b_{{1}}^{2}+50\,b_{{3}}^{2}+26\,b_{{1}}b_{{3}}
+ \mathcal{O}(b_1a_0, b_3a_0; \, a_0, b_1, b_3).
\end{array}\eeq
It follows that $\, \delta \phi \equiv 0, \mod \; a_0, \, b_1, \, b_3,$
and that the structure equations for $ \{ \, a_0, \, b_1, \, b_3 \, \}$ are closed, i.e.,
their derivatives are expressed as the functions of themselves
and do not involve any new variables.

\begin{remark}\label{3rmk}
Note that
$$da_0 \equiv -4 (b_1 \omega^1 - b_3 \omega^2), \mod \; \theta; \,a_0.$$
Hence if $\, a_0$ vanishes up to order one at a point on a connected open subset $\, M$,
the uniqueness theorem of ODE implies that $\, a_0, \, b_1, \, b_3 \equiv 0$ identically.
Hence in this case $\, \delta \phi = 0$ and the deformation is trivial.
\end{remark}
\begin{proposition}\label{3propo}
Let $\, \W$ be a linear Legendrian 3-web on a connected open subset $\, M \subset \PP^3$.
Consider a linear Legendrian deformation of $\, \W$ represented by
the $\, \mathfrak{sp}_2\C$-valued 1-form $\, \delta \phi$
satisfying the initial conditions \eqref{3step00}, \eqref{3step01}.
Then the components of $\, \delta \phi$ are given by,
and satisfy the structure equations \eqref{3step0set}, \eqref{3step0set0},
\eqref{3step0set1}, \eqref{3step0set2}, \eqref{3step0set3}.

Suppose $\, \delta \phi \equiv 0,\mod \; \theta$ at a point in $\, M$.
Then $\, \delta \phi$ vanishes identically and the deformation is trivial.
\end{proposition}
\Pf If $\, \delta \phi \equiv 0,\mod \, \theta\,$ at a  point,
\eqref{3step0set} shows that $\, a_0, \, b_1, \, b_3=0$ at the given point.
The rest follows from Remark \ref{3rmk}.
$\sq$

\subsection{Bound on the number of distinct  linearizations}\label{sec31}
The remaining compatibility conditions from
$\,\Delta^0_1, \, \Delta^0_3, \Delta^0_2; \,d(d(b_1))=0,\, d(d(b_3))=0$,
impose a set of polynomial equations (and their successive derivatives)
on the deformation parameters $\, \{\,a_0,\,b_1,\,b_3\,\}$.
The analysis of the root structure of these equations inevitably leads to
a variety of case by case analysis problems depending on
the relative invariants of the original 3-web $\, \W$,
e.g., the resultants of a set of polynomial compatibility equations for
$\, \{\,a_0,\,b_1,\,b_3\,\}$ are expressed in terms of
the local invariants of $\, \W$.

Due to the complexity and size of the algebraic analysis involved,
we shall consider the first few lowest order compatibility equations.
In this section, we examine these equations without any extra conditions on the local
invariants of the original web $\,\W$,
and determine an upper bound on the number of distinct linearizations.

\two
We continue the analysis from \textbf{Step 3'}.

The identities $\,d(d(b_1))\w\theta=0,\, d(d(b_3))\w\theta=0$ give
a set of two compatibility equations which must vanish identically
for a linear Legendrian deformation.
\begin{align}\label{3Eq13}
Eq_1&\equiv b_1^2+2b_1b_3,  \\
Eq_3&\equiv b_3^2+2b_1b_3, \mod\;
\textnormal{\{terms at most linear in $\, b_1, b_3$\}.} \n
\end{align}
$\, Eq_1, \, Eq_3$ are polynomials of degree 4 \,in $\, \{\,a_0,\,b_1,\,b_3\,\}$.
One computes then that  $\, \Delta^0_1, \, \Delta^0_3 \equiv 0\,$ modulo $\, Eq_1, \, Eq_3$.

$\Delta^0_2$ modulo $\, Eq_1, \, Eq_3$ gives another set of two equations
\begin{align}\label{3Eq0}
Eq_0 &\equiv  a_0^5+\frac{314}{111} A_0 a_0^4, \\
Eq_9 &\equiv  A_0 a_0^4,\quad\quad\qquad\mod\;
\textnormal{\{terms at most linear in $\, b_1, b_3$, and of degree $\leq 3$ in $\, a_0$ \}.}\n
\end{align}
$Eq_0, \, Eq_9$ are polynomials of degree 5, 4\, in $\, \{\,a_0,\,b_1,\,b_3\,\}$
respectively.

\begin{theorem}\label{31thm}
Let $\, \W$ be a Legendrian 3-web on a connected contact three manifold $\, M$.
Then it admits at most $\,4 \cdot 4\cdot 5 +1=81$ distinct local linearizations.
If the local invariant $\, A_0$ of a linearization of $\, \W$ is nonzero,
it admits at most $\,4 \cdot 4\cdot 4 +1=65\,$ distinct local linearizations.\footnotemark
\footnotetext{In case the invariant $\, A_0 \equiv 0$ identically,
the linear Legendrian 3-web admits at most two distinct local linearizations,
see Section \ref{sec4}.}
\end{theorem}
\Pf
It is clear that $\, Eq_1$ and $\, Eq_3$ cannot have a common linear factor.
The rest follows by counting the degrees of the polynomials $\,\{ Eq_1, \, Eq_3,\,Eq_0\}$,
and $\, \{ Eq_1, \, Eq_3,\,Eq_9 \}$.
$\sq$

It is unlikely that this is the optimal bound.
We suspect that the optimal bound on the number of distinct local linearizations is two,
see Conjecture \ref{conj} in Section \ref{sec5}.
Theorem \ref{31thm} shows that the number of distinct local linearizations
of a Legendrian 3-web is uniformly bounded.

\subsection{Rigidity of the linear Legendrian 3-webs with a flat point}\label{sec32}
Remark \ref{3rmk} implies that a linear Legendrian 3-web is rigid
under linear Legendrian deformation
when the various compatibility equations force the deformation parameters
$\, \{a_0,\,b_1,\,b_3\}$ to vanish at a single point.
In this section we examine $\,\{ Eq_1, \, Eq_3,\,Eq_0\}$ in \eqref{3Eq13}, \eqref{3Eq0},
and show that this occurs in case $\, \W$ is sufficiently flat at a point.

\one
Fix a reference point $\, \textnormal{x}_0 \in M$.
Assume all the coefficients $\, A_i, \, B_j, \, C_k$'s and
their derivatives of sufficiently high order
vanish at $\, \textnormal{x}_0$.
The compatibility equations $\, Eq_1, \, Eq_3, \, Eq_0$ in \eqref{3Eq13}, \eqref{3Eq0}
evaluated at $\, \textnormal{x}_0$ become (up to scaling by constants)
\beq\label{5reduced}
\begin{array}{rl}
(Eq_1-Eq_3)|_{\textnormal{x}_0}&=(b_1-b_3)(4b_1+4b_3-11 a_0^2),   \\
Eq_1|_{\textnormal{x}_0}&=(4 b_1+5 a_0^2)(2 b_1+4 b_3-3 a_0^2),   \\
Eq_0|_{\textnormal{x}_0}&=a_{{0}} \left(185 a_0^4-60 a_0^2 b_3-220 b_1a_0^2
-208 b_1b_3 \right).
\end{array}
\eeq
It is easily checked that the only root to this system of equations
for $\, \{\,a_0, \,b_1, \, b_3\,\}$ is $\; a_0,\,b_1,\,b_3|_{\textnormal{p}_0} = 0$.
By Remark \ref{3rmk}, in this case $\, \{\,a_0, \,b_1, \, b_3\,\}$ vanish identically
and the deformation is trivial.

The following theorem describes up to which order the coefficients
$\, A_i, \, B_j, \, C_k$'s should vanish at the reference point
to put the equations $\, Eq_1, \, Eq_3, \, Eq_0$ into the reduced form \eqref{5reduced}.
\begin{theorem}\label{32thm}
Let $\, \W$ be a linear Legendrian 3-web
on a connected open subset $\, M \subset \PP^3$.
Let $\, A_i, \, B_j, \, C_k$ be the structure coefficients of $\, \W$, \eqref{2phi}.
Let $\, \textnormal{x}_0 \in M$ be a reference point.
Suppose
$$\begin{array}{rl}
A_0& \textnormal{vanishes to order 3 at p$_0$}, \n \\
A_1, \, A_2, \,B_1,\,B_2,\,C_1,\,C_2&
\textnormal{vanish  to order 1 at p$_0$}, \n \\
B_0, \, C_9& \textnormal{vanish at p$_0$}. \n
\end{array}$$
Then $\, \W$ is rigid and does not admit any nontrivial
local linear Legendrian deformations.
\end{theorem}
\Pf
Examining the equations $\, Eq_1, \, Eq_3, \, Eq_0$, the vanishing conditions
on $\, A_i, \, B_j, \, C_k$ are sufficient to imply \eqref{5reduced}.
The rest follows from Remark \ref{3rmk}.
$\sq$

This agrees with the partial proof of the Gronwall conjecture for the planar 3-webs
obtained in \cite{Wa}.

\section{Linearization of the Legendrian 3-webs of maximum rank}\label{sec4}
Theorem \ref{32thm} states that a linear Legendrian 3-web is rigid
when it is  sufficiently flat at a point.
A question arises as to which abstract Legendrian 3-webs
are such that their linearizations are likely to have the similar property.
The first candidates would be the Legendrian 3-webs of maximum rank three, \cite{Wa2}.
The following refinement of Theorem \ref{32thm} gives
a partial proof of the Legendrian Gronwall Conjecture, Section \ref{sec5},
for this class of Legendrian 3-webs.

\begin{theorem}\label{4thm}
Let $\, \W$ be a Legendrian 3-web on a connected contact three manifold $\, M$.
Suppose $\, \W$ has the maximum rank three.

a)
$\, \W$ admits a local linearization $\, M \hook \PP^3$
as the dual Legendrian 3-web of an analytic surface
$\, \Sigma \subset \Q^3$
which is the union of three hyperplane sections
$\, \Sigma=\cup_{i=1}^3 H^i, \; H^i\subset\Q^2$.
Conversely, for any set of three distinct hyperplane sections $\,H^i$ in $\, \Q^3$,
the dual Legendrian 3-web $\,\W_{\cup_{i=1}^3H^i}$ has the maximum rank.

b)
The linearization of $\, \W$  in a) is unique up to motion by \,Sp$_2\C$,
with the only exception when the structure invariants of $\, \W$ in \eqref{webstrt}
satisfy the relation
\beq\label{4STR}
 S = -T = 2R \ne 0.
\eeq
In this case, the Legendrian 3-web admits
exactly one more distinct local linearization
as the dual Legendrian 3-web of the Legendrian twisted cubic.
\end{theorem}
\begin{corollary}
A Legendrian 3-web of maximum rank is algebraic.
\end{corollary}

\emph{Proof \,of \,Theorem \ref{4thm}.}

a) Let $\, \W$ be a Legendrian 3-web of maximum rank three.
Then the structure invariants of $\, \W$ in \eqref{webstrt} satisfy the relations
\begin{align}
dR,\, \, dS, \, dT &\equiv 0, \mod \; \rho, \n \\
L = K &= 0.\n
\end{align}
Recall \eqref{2linearQ} in Section \ref{sec2}
for the linear Legendrian 3-webs with the vanishing local invariant $\,A_0$.
Substituting
\beq\label{4phi}
 B_0 = \frac{R}{4},\;A_1=-\frac{T}{2},\;A_2=-\frac{S}{2};\; \, \rho_0 = -\rho,
\eeq
one has $\, d\phi + \phi \w \phi=0$,
and it induces a local linearization of $\, \W$.

\one
b) Given the linearization of $\, \W$ defined by $\, \phi$ in \eqref{2linearQ}, \eqref{4phi},
the deformation analysis as in Section \ref{sec3} shows that
the various compatibility equations force $\, \delta \phi =0$,
except when $\, S = - T = 2R \ne 0$
(this requires a long but straightforward case by case analysis).
In this case $\, \W$ admits exactly one more distinct linearization
induced by the (deformed) Maurer-Cartan form
\beq\label{4phi'}
\phi'=
\left[ \begin {array}{cccc} -\rho&\frac{3}{2} R \omega^1+\frac{3}{4} R \omega^2
&-{\frac {27}{8}} R^{2}\theta&\frac{3}{4} R\omega^1+\frac{3}{2} R \omega^2
\\\noalign{\medskip}\omega^1&\frac{1}{2}  R \theta+a_{{0}}\omega^1+a_{{0}}\omega^2
&\frac{3}{4} R \omega^1+\frac{3}{2} R \omega^2&R \theta+a_{{0}}\omega^1
\\\noalign{\medskip}2\,\theta&\omega^2&\rho&-\omega^1
\\\noalign{\medskip}\omega^2&-R \theta-a_{{0}}\omega^2&-\frac{3}{2} R \omega^1
-\frac{3}{2} R \omega^2&-\frac{1}{2} R \theta-a_{{0}} \omega^1-a_{{0}}\omega^2\end {array}
 \right],
\eeq
where $\, a_0^2 = 2 R$.

We claim that this describes the dual 3-web of the Legendrian twisted cubic.
Choose the linear Legendrian foliation $\, \mathcal{F}^1$ defined by
$\, \langle \, \omega^2, \, \theta \, \rangle^{\perp}$ via $\, \phi'$.
The corresponding surface $\, \Sigma^1 \subset \Q^3$
is described by $\, [ Z_0 \w Z_1 ]$.
The conformal structure of $\,\Q^3$ is represented by the quadratic form
$\,\beta^0_0\beta^1_1-(\beta^1_0)^2$, whose restriction on $\,\Sigma^1$
becomes the perfect square $\,(\omega^2-a_0\theta)^2$.
Differentiating $\,Z_0 \w Z_1$, one gets that
the two dimensional subspace spanned by
$\,\{Z_0\w Z_1,\,(2Z_1-a_0Z_0)\w(Z_2+\frac{a_0}{2}Z_3)\}$ is constant
along the leaves of the foliation on $\,\Sigma^1$ defined by $\,(\omega^2-a_0\theta)^{\perp}$.
By duality, $\, \Sigma^1$ is ruled by the null line dual to
$\, \textnormal{p}_1=[ 2Z_1-a_0Z_0 ]\in\PP^3.$
The similar analysis for the foliations
$\, \mathcal{F}^2, \, \mathcal{F}^3$ shows that the corresponding surfaces
$\,\Sigma^2,\,\Sigma^3$ are ruled by the null lines dual to
$\, \textnormal{p}_2=[ 2Z_3+a_0Z_0 ],\,\textnormal{p}_3=[ 2(Z_1-Z_3)+a_0Z_0 ]$
respectively. Note that $\,\cap_{i}\PP^2_{\textnormal{p}_i} = [Z_0]$.

Our claim is that the loci of $\,\textnormal{p}_1,\,\textnormal{p}_2,\,\textnormal{p}_3$,
all lie in the same Legendrian twisted cubic. A short computation shows that
each $\,\textnormal{p}_i$ describes a Legendrian curve.
One may then verify by direct computation that the following three quadratic polynomials
generate the well defined, covariant constant,
three dimensional space of quadratic polynomials in
a neighborhood of $\,[Z_0]$ which vanish simultaneously on $\,\textnormal{p}_i$'s.
Here $\,W=(W_0,\,W_1,\,W_2,\,W_3)^t$ is the dual frame of $\,Z=(Z_0,Z_1,Z_2,Z_3)$
that satisfies the structure equation $\,dW=-\phi'W$.
$$
\begin{array}{rl}
Q_0&=-R  W_{{1}}^{2} - R W_{{1}}W_{{3}} +2\, W_{{0}}^{2}-R W_{{3}}^{2}
+{\frac {27}{8}}R^2 W_{{2}}^{2} , \n \\
Q_1&=a_{{0}} W_{{1}}^{2}+2\,W_{{0}}W_{{1}}-3R W_{{2}}W_{{3}}
+2a_{{0}} W_{{1}}W_{{3}} -\frac{3}{2} R W_{{1}}W_{{2}}, \n \\
Q_3&=-2a_{{0}}W_{{1}}W_{{3}} -a_{{0}} W_{{3}}^{2} +2\,W_{{0}}W_{{3}}
+\frac{3}{2} R W_{{2}}W_{{3}} +3R W_{{1}}W_{{2}}. \n
\end{array}
$$
It is clear that these polynomials define the Legendrian twisted cubic in $\,\PP^3$
that connects the three arcs traversed by $\,\textnormal{p}_i$'s.
$\sq$

\two
The analysis in the proof of b) above has the following algebro-geometric implication.
Let $\,\W$ be the dual 3-web of the Legendrian twisted cubic just described.
Since the Gauss map of the Legendrian twisted cubic in $\, \PP^3$
is the null rational normal curve in $\, \Q^3 \subset \PP^4$,
$\, \W$ is also the dual 3-web of the analytic surface  $\, \Sigma \subset \Q^3$
which is the tangent developable of the null rational normal curve.

Consider the following analogue of the converse of Abel's theorem, \cite{CG1,HP}.

\vspace{4mm}
\parbox{16.5cm}{
\fbox{Legendrian analogue of the converse of Abel's theorem}

Let $\,\x_0\in\PP^3$ be a reference point. Let $\,N_0\subset\Q^3$ be the null line
dual to $\,\x_0$. Let $\,\Sigma^i\subset\Q^3,\,i=1,\,2,\,...\,d,$ be a set of
$\,d$ distinct pieces of local analytic surfaces each of which intersects
$\,N_0$ transversally at a single point $\,\textnormal{q}_i(\x_0)$.
Let $\,\Omega_i$ be a meromorphic 1-form on $\,\Sigma^i$ which is regular
at $\,\textnormal{q}_i(\x_0)$. Suppose the following local trace vanishes
in a neighborhood of $\,\x_0$.
$$
\textnormal{Tr}\,(\,\Omega_i\,)=\sum_i \textnormal{q}_i^*\Omega^i =0.
$$
Then there exists a null degree $\,d$ analytic surface $\,\Sigma$,
and a meromorphic 1-form $\,\Omega$ on $\,\Sigma$, which analytically extends the given
data $\,\cup_i(\,\Sigma^i,\,\Omega_i)$.
}

\vspace{4mm}
Assuming this is true, the relevant observation is that
\emph{
the tangent developable of the null rational normal curve
lies in the intersection of $\,\Q^3$ with a cubic hypersurface of $\,\PP^4$.}\footnotemark
\footnotetext{Let $\,V_m$ be the irreducible SL$_2\C$-module of dimension $\,m+1$.
By Clebsch-Gordan, the symmetric cubic tensor product decomposes into
$\,S^3(V_4)=V_{12}\oplus V_8\oplus V_6\oplus V_4\oplus V_0$.
The $\,V_0$ piece vanishes on the tangent developable of
the null rational normal curve.} It is known that a smooth complete intersection
of type $\, (2, \, 3)$ in $\, \PP^4$ is a K3 surface, which has no nonzero holomorphic 1-forms.
According to Theorem \ref{4thm},
when this complete intersection degenerates to $\, \Sigma$,
it supports three closed generalized holomorphic 1-forms.
We currently do not have any purely algebro-geometric explanation of this phenomenon.

\section{Concluding remarks}\label{sec5}
1.
The Legendrian analogue of the Gronwall conjecture can be stated as follows.

\begin{conjecture}[Legendrian Gronwall conjecture]\label{conj}
Let $\, \W$ be a Legendrian 3-web on a connected contact three manifold.
Then it admits at most one distinct local linearization in $\, \PP^3$,
with the only exception when $\, \W$ is locally equivalent to
the dual 3-web of the Legendrian twisted cubic curve in $\, \PP^3$
(in which case $\, \W$ admits exactly two distinct linearizations).
\end{conjecture}
Note the discrepancy when compared with the planar web case.
An algebraic planar 3-web admits infinitely many distinct local linearizations.

From the analysis carried out in this paper,
one suspects that a generic Legendrian 3-web admits at most one distinct linearization,
although we did not write down the condition for rigidity explicitly.
Since the other extreme case of the maximum rank Legendrian 3-webs are treated
in Section \ref{sec4},
the next step toward the proof of the Legendrian Gronwall conjecture
would be to consider the class of Legendrian 3-webs
with one, or two independent Abelian relations.
One may then attempt to show that a Legendrian 3-web with two distinct linearizations
necessarily possesses at least one Abelian relation.
It appears to be a difficult problem to directly analyze the root structure of
the polynomial integrability equations for the linear Legendrian deformation
discussed in Section \ref{sec3}.

\two
2.
Let $\, \Sigma \subset \Q^3$ be a null degree $\, d$ surface.
We say that $\, \Sigma$ is \emph{extremal}
when the associated dual linear Legendrian $\, d$-web $\, \W_{\Sigma}$
has the maximum rank $\, \rho_d=\frac{(d-1)(d-2)(2d+3)}{6}$, \cite{Wa2}.
Theorem \ref{4thm} implies that
the only null degree 3 extremal surfaces are
the union of three hyperplane sections, and the tangent developable to
the null rational normal curve.
Can one give an algebro-geometric proof of this?
The direct differential analysis for the Legendrian $\, d$-webs of maximum rank
for the case $\, d \geq 4$ is complicated.
One may hope to generalize the algebro-geometric classification of
the extremal null degree 3 surfaces to the extremal null degree $\, d$ surfaces
in general.

One may start by considering the following sub-problem.
Let $\,\gamma\subset\PP^3$ be a degree $\,d$ curve.
We say that $\, \gamma$ is \emph{extremal}
when the associated dual linear Legendrian $\, d$-web $\, \W_{\gamma}$
has the maximum rank $\, \rho_d$.
Theorem \ref{4thm} shows that
it is possible for $\, \W_{\gamma}$ to have an Abelian relation
which is not induced from a holomorphic 1-form on $\, \gamma$.
We suspect that this kind of auxiliary Abelian relations exist
only when $\, \gamma$ is itself Legendrian.
Can one give an algebro-geometric classification of the extremal Legendrian curves?
As mentioned in Section \ref{sec1}, such an extremal Legendrian curve
is tightly controlled by a large number of generalized addition laws.

\two
3.
As discussed briefly in \cite{Wa2},
the geometry of a Legendrian 2-web is locally equivalent to
the geometry of a single scalar second order ODE up to point transformation.
In particular, the structure of a Legendrian 2-web has local invariants,
and not every two of them are locally equivalent, \cite{Ca}.

A question arises as to the geometric meaning of the linearization
of a Legendrian 2-web.
Can the linearizability be considered as a local counterpart of the notion of completeness
of the associated projective connection discussed in \cite{McK}?

\renewcommand{\theequation}{A-\arabic{equation}}
\setcounter{equation}{0}  
\section*{Appendix}
Let $\, \W$ be a Legendrian 3-web on a contact three manifold $\, M$.
There exists a sub-bundle $\, B$ of the $\,  \textnormal{GL}_3 \C$
principal frame bundle  of $\, M$ on which the tautological 1-forms
$\, \theta; \, \omega^i, \, i = 1, \, 2, \, 3, \,\sum\omega^i=0$,
satisfy the following structure equations.
\begin{align}\label{webstrt}
d
\begin{pmatrix}
\omega^1 \\ \omega^2 \\  \theta
\end{pmatrix}
&= -
\begin{pmatrix}
\rho & \cdot &  \cdot \\
\cdot & \rho & \cdot \\
\cdot &  \cdot &  2 \rho
\end{pmatrix}
\w
\begin{pmatrix}
\omega^1 \\ \omega^2 \\  \theta
\end{pmatrix}
+
\begin{pmatrix}
\theta \w ( R \, \omega^1 + S \, \omega^2 )  \\
\theta \w ( T \, \omega^1 - R \, \omega^2 )  \\
\omega^1 \w \omega^2
\end{pmatrix},\\
d \rho &= \theta \w ( \,L \omega^1 + K \, \omega^2 \, ).\n
\end{align}
Here $\, R, \, S, \, T,\,L,\,K$ are torsion coefficients.
The contact structure on $\,M$ is defined by $\, \langle\,\theta\,\rangle^{\perp}$.
$\,\W$ is defined by the three  line fields $\,\langle\,\omega^i, \, \theta\,\rangle^{\perp}$.
See \cite{Wa2} for the derivation of this structure equation.


\end{document}